\documentclass[12pt]{amsart}

\newtheorem{theorem}{Theorem}[section]
\newtheorem{lemma}[theorem]{Lemma}

\theoremstyle{definition}

\theoremstyle{remark}

\numberwithin{equation}{section}

\usepackage{times}
\usepackage{enumerate}
\usepackage{tfrupee}
\usepackage{tikz}
\usetikzlibrary{chains,fit}
\usetikzlibrary{shapes,snakes}
\usetikzlibrary{graphs}
\usepackage{graphicx,adjustbox}
\usepackage{hyperref}
\usepackage{amsmath,amssymb}

\usepackage{amscd}
\usepackage{graphicx}
\usepackage[all]{xy}

\title[Quadrics defined by skew-symmetric matrices]
{Quadrics defined by skew-symmetric matrices}
\author{
Joydip Saha
\and
Indranath Sengupta
\and
Gaurab Tripathi
}
\date{}

\address{\small \rm  Discipline of Mathematics, IIT Gandhinagar, Palaj, Gandhinagar, 
Gujarat 382355, INDIA.}
\email{saha.joydip56@gmail.com}
\thanks{The first author is a research associate supported by the 
research project EMR/2015/000776 sponsored by the SERB, Government of India.}

\address{\small \rm  Discipline of Mathematics, IIT Gandhinagar, Palaj, Gandhinagar, 
Gujarat 382355, INDIA.}
\email{indranathsg@iitgn.ac.in}
\thanks{The second author is the corresponding author, who is supported by the 
research project EMR/2015/000776 sponsored by the SERB, Government of India.}

\address{\small \rm Department of Mathematics, Jadavpur University, Kolkata,
WB 700 032, India.} 
\email{gelatinx@gmail.com}
\thanks{The third author thanks CSIR for the Senior Research Fellowship.}

\date{}

\subjclass[2010]{Primary 13D02; Secondary 13C40, 13P10, 13D07.}

\keywords{Gr\"{o}bner basis, Betti numbers, determinantal ideals, skew-symmetric matrix, 
mapping cone.}

\begin{document}
\begin{abstract}
In this paper we propose a model for computing a minimal free resolution for 
ideals of the form $I_{1}(X_{n}Y_{n})$, where $X_{n}$ is an $n\times n$ 
skew-symmetric matrix with indeterminate entries $x_{ij}$ and $Y_{n}$ is a generic 
column matrix with indeterminate entries $y_{j}$. We verify that the model works 
for $n=3$ and $n=4$ and pose some statements as conjectures. Answering the conjectures 
in affirmative would enable us to compute a minimal free resolution for general $n$.
\end{abstract}

\maketitle

\section{Introduction}
Let $K$ be a field. Let $\{x_{ij}; \, 1\leq i \leq m, \, 1\leq j \leq n\}$, 
$\{y_{j}; \, 1\leq j \leq n\}$ be indeterminates over $K$, so that 
$R=K[x_{ij}, y_{j}]$ denotes the polynomial algebra over $K$. Let $X_{n}$ 
denote an $n\times n$ skew-symmetric matrix such that its entries are the indeterminates 
$\pm x_{ij}$ and $0$. We call such a matrix a generic skew-symmetric matrix. 
Let $Y_{n}=(y_{j})_{n\times 1}$ be the generic $n\times 1$ matrix. 
It is very hard to compute a graded minimal free resolution of the ideal 
$I_{1}(X_{n}Y_{n})$. 
\medskip

Ideals of the form $I_{1}(X_{n}Y_{n})$ has been studied by \cite{herzog} and 
they appear in some of our recent works; see \cite{sstgbasis}, \cite{sstprimary}, \cite{sstsum}, 
\cite{sstregseq}. We described its Gr\"{o}bner bases, primary decompositions 
and Betti numbers through computational techniques, mostly under the assumption 
that $X_{n}$ is either a generic or a generic symmetric matrix. It is in deed 
the case that these ideals are far more difficult to understand when $X_{n}$ 
is a generic skew-symmetric matrix. 
\medskip

In this paper, we present a scheme for computing a graded minimal free resolution 
of the ideal $I_{1}(X_{n}Y_{n})$, where 
$X_{n}$ is a $n\times n$ generic skew-symmetric matrix and and $Y_{n}$ is a 
generic $n\times 1$ matrix. We show that if we assume the truth of two statements 
then the scheme works for a general $n$. These two statements which have been proposed 
as conjectures appear to be correct as seen from symbolic computation using the 
computer algebra software \textit{Singular} \cite{DGPS}. We finally verify the validity 
of these conjectures for $n=3$ and $n=4$. We refer to \cite{peeva} for basic knowledge 
on the techniques used by us. 

\section{General Scheme and Conjectures}

Let $X_{n}=\left[
\begin{array}{ccccc}
0& x_{12} &x_{13} & \ldots &x_{1n}\\
-x_{12}&0 & x_{23}& \ldots &x_{2n}\\
x_{13}& x_{23}&&\\
\vdots& \vdots &&\\
-x_{1n}& -x_{2n}&\ldots && 0
\end{array}
\right]
$ and $Y_{n}
 =\left[
\begin{array}{c}
y_{1}\\
\vdots\\
y_{n}
\end{array}
\right]
$. 
\medskip

\noindent Our aim is to find a minimal free resolution of $I_{1}(X_{n}Y_{n})$. 
Assuming $x_{ij}=-x_{ji}$, if $i>j$ and $x_{ii}=0$, let $g_{ki}=\Sigma_{j=1}^{i}x_{kj}y_{j}$. 
Therefore the generators of $I_{1}(X_{n}Y_{n})$ are 
$\langle g_{1i},g_{2i},\cdots, g_{ni}\rangle$.  Let $\Delta_{(i)n}$ 
denote the Pffafian of the skew symmetric matrix $X_{n}$ with 
the $i$-th row and the $i$-th column deleted.
\medskip

\begin{lemma}\label{useful} 
Assuming $x_{ij}=-x_{ji}$, if $i>j$ and $x_{ii}=0$
\begin{itemize}
\item[(i)] $y_{n}g_{nn} = -\left( y_{1}g_{1n}+ y_{2}g_{2n}+\cdots+ y_{n-1}g_{(n-1)n}\right)$.

\item[(ii)] $g_{k (n-1)}g_{nn}  = x_{kn}y_{1}g_{1n}+ x_{kn}y_{2}g_{2n}+ \cdots + g_{nn}+x_{kn}y_{k}g_{kn}+\cdots + x_{kn}y_{n}g_{n-1 n}$.

\item[(iii)] $ \Delta_{(n)n}y_{n}  =  (-\Delta_{(1)n})g_{1n}+ (\Delta_{(2)n})g_{2n}+ \cdots + ((-1)^{n-1}\Delta_{(n-1)n})g_{(n-1)n}$.
\end{itemize}
\end{lemma}

\proof A simple calculation gives the proof. \qed
\medskip

\begin{lemma}\label{akosul}
$\{g_{1n},g_{2n},\cdots , g_{(n-1) n} \}$ forms a regular sequence for $n\geq 2$.
\end{lemma}

\proof See part (ii) of Theorem 2.2 in \cite{sstprimary}. \qed
\bigskip

\noindent\textbf{Notations.}
\begin{enumerate}
\item[(i)]Let $I_{n}=\langle g_{1n},g_{2n},\cdots , g_{(n-1) n}\rangle$. 
By Lemma \ref{akosul}, the ideal $I_{n}$ is minimally 
resolved by the Koszul complex
$$0\longrightarrow R \longrightarrow \cdots \longrightarrow R^{\binom{n-1}{2}}\stackrel{\psi_{2n}}{\longrightarrow}R^{\binom{n-1}{1}}\stackrel{\psi_{1n}}{\longrightarrow}R \stackrel{\psi_{0n}}{\longrightarrow} R/I_{n}\longrightarrow 0 ;$$
where $\psi_{kn}: R^{\binom{n-1}{k}}\longrightarrow R^{\binom{n-1}{k-1}}$ and  $k\in\{0,1, 2,\cdots, n-1 \}$.

\item[(ii)] Let $J_{n}= \langle g_{nn}\rangle$ and $L_{n}= I_{n}+J_{n}=I_{1}(X_{n}Y_{n})$.
\end{enumerate}
\medskip

\noindent Computations with Singular give us enough evidence in support of the Conjectures 
proposed below:
\medskip

\noindent\textbf{Conjecture 1.} $C_{n}:= (I_{n}:J_{n})=\langle g_{1 (n-1)}, g_{2 (n-2)},\cdots, g_{n-1 (n-1)}, y_{n}, \Delta_{(n)n}\rangle $. If $n$ is even then $\Delta_{(n)n}=0$ and $C_{n}=\langle g_{1 (n-1)}, g_{2 (n-2)},\cdots, g_{n-1 (n-1)}, y_{n}\rangle$, for every $n\geq 4$.
\medskip

\noindent\textbf{Conjecture 2.} If $n$ is odd then $\Delta_{(n)n}\neq 0$. For every $n\geq 4$,
$$P_{n}:= (\langle g_{1 (n-1)}, g_{2 (n-2)},\cdots, g_{n-1 (n-1)}\rangle:\Delta_{(n)n})
= (y_{1},\cdots, y_{n-1}).$$
\medskip

\noindent Assuming the validity of these conjectures we can construct a minimal free 
resolution for $I_{1}(X_{n}Y_{n})$ through the following steps.
\medskip

We proceed by induction on $n\geq 3$. We first compute a resolution of $L_{3}$, which is 
not difficult. For 
$3\leq i-1<n$, let a resolution of $L_{i-1}$ be 
$$0 \longrightarrow \cdots \longrightarrow R^{\beta_{1(i)}}\stackrel{d_{1i}}{\longrightarrow}\longrightarrow R^{\beta_{0(i)}}\stackrel{d_{0i}}{\longrightarrow} R/L_{i-1}\longrightarrow 0 $$ where $d_{0i}: R\longrightarrow R/L_{i-1}$ is the projection map.
\medskip

A resolution of $P_{n}$ is the Koszul complex, which is the following:
$$0\longrightarrow R \longrightarrow \cdots \longrightarrow R^{\binom{n-1}{2}}\stackrel{\phi_{2n}}{\longrightarrow}R^{\binom{n-1}{1}}\stackrel{\phi_{1n}}{\longrightarrow}R \stackrel{\phi_{0n}}{\longrightarrow} R/P_{n}\longrightarrow 0 $$
where $\phi_{kn}: R^{\binom{n-1}{k}}\longrightarrow R^{\binom{n-1}{k-1}}$, for 
$k\in\{0,1, 2,\cdots, n-1 \}$.
\medskip

\noindent\textbf{Case 1.} For $i< n$ and $i$ is odd, let 
$T_{i}:= L_{i-1}+ \langle\Delta_{(i)i}\rangle$.  Using mapping cone we get,
\medskip

\xymatrix{
\cdots \ar[r] &R^{\binom{i-1}{3}}\ar[d]^{\delta_{3i}}\ar[r]& R^{\binom{i-1}{2}}\ar[r]^{\phi_{2i}}\ar[d]^{\delta_{2i}}&
R^{\binom{i-1}{1}}\ar[r]^{\phi_{1i}}\ar[d]^{\delta_{1i}}&R\ar[r]^{\phi_{0i}}\ar[d]^{\delta_{0i}= \Delta_{i(i)}} &R/P_{i}\ar[d]\ar[r]&0 \\
\cdots \ar[r] &R^{\beta_{3(i-1)}}\ar[r] &R^{\beta_{2 (i-1)}}\ar[r]_{d_{2}i-2}&R^{\beta_{1(i-1)}}\ar[r]_{d_{1}i-1}&R\ar[r]_{d_{0i}}&R/L_{i-1}\ar[r]& 0
}
\medskip

\noindent Therefore a resolution of $T_{i}$ is 
$$\longrightarrow R^{\binom{i-1}{1}}\oplus R^{\beta_{2(i-1)}}\stackrel{\delta_{2i}^{'}}{\longrightarrow} R\oplus R^{\beta_{1(i-1)}}\stackrel{\delta_{1i}^{'}}{\longrightarrow} R\longrightarrow R/T_{i}\longrightarrow 0; $$
where $\delta_{ki}^{'}: R^{\binom{i-1}{k-1}} \oplus R^{\beta_{k(k\i-1)}}\longrightarrow  R^{\binom{i-1}{k-2}} \oplus R^{\beta_{k(k\i-1)}} $  and $\delta_{ki}^{'} = \left[
\begin{array}{cc}
-\phi_{k-1 i} & 0 \\
\delta_{k-1 i} & d_{k i-1}
\end{array}
\right]$.
\medskip

The resolution of $T_{i}$ obtained above may not be minimal. 
Assuming that we can extract a minimal free resolution from 
it by identifying matching of degrees and cancelling them 
(see the computations for some special values of $n$ in the 
next section) let $\mathbb{T}_{i}$ be the minimal free 
resolution of $T_{i}$, whose differentials are $\delta_{ki}$, i.e. 
$$\mathbb{T}_{i}\cdots \longrightarrow R^{\gamma_{2i}}\stackrel{\delta_{2i}}{\longrightarrow}R^{\gamma_{1i}}\stackrel{\delta_{1i}}{\longrightarrow} R\longrightarrow R/T_{i}\longrightarrow 0$$
\medskip
 
To find the resolution of $C_{i}$, we need to tensor 
$0\longrightarrow R\stackrel{y_{i}} {\longrightarrow} R\longrightarrow 0$ 
with the complex $\mathbb{T}_{i}$, which gives us

$\cdots\longrightarrow (R^{\gamma_{3i}}\otimes R)\oplus (R^{\gamma_{2i}}\otimes R)\stackrel{\eta_{3i}}{\longrightarrow} (R^{\gamma_{2i}}\otimes R)\oplus (R^{\gamma_{1i}}\otimes R)\stackrel{\eta_{2i}}{\longrightarrow} (R^{\gamma_{1i}}\otimes R)\oplus (R\otimes R)\stackrel{\eta_{1i}}{\longrightarrow} 
R\otimes R \longrightarrow R/C_{i}\longrightarrow 0$;

\noindent where $\eta_{ki}= \left[
\begin{array}{cc}
\delta_{ki} & -y_{i}I \\
0 & \delta_{k-1 i}
\end{array}
\right]$.\hfill(*)
\medskip

We first rewrite complex(*), which gives us a minimal free resolution of $C_{i}$. 
Then, we construct the mapping cone of the following complexes 
with respect to the following connecting maps:
\medskip

\xymatrix{
\cdots \ar[r] &R^{(\gamma_{i3}+\gamma_{i2})}\ar[d]^{\xi_{3i}}\ar[r]& R^{(\gamma_{i2}+\gamma_{i1})}\ar[r]^{\eta_{2i}}\ar[d]^{\xi_{2i}}&
R^{(\gamma_{i1}+1)}\ar[r]^{\eta_{1i}}\ar[d]^{\xi_{1i}}&R\ar[r]^{\eta_{0i}}\ar[d]^{\xi_{0i}= g_{(i)i}} &R/C_{i}\ar[d]\ar[r]&0 \\
\cdots \ar[r] &R^{\binom{i-1}{3}}\ar[r] &R^{\binom{i-1}{2}}\ar[r]_{\psi_{2(i)}}&R^{\binom{i-1}{1}}\ar[r]_{\psi_{1(i)}}&R\ar[r]_{\psi_{0i}}&R/I_{i}\ar[r]& 0
} 

\noindent We create a minimal resolution out after the mapping cone construction by 
suitable cancellation of matched degrees. 
\medskip

\noindent\textbf{Case 2.} Let $i$ be even and $i< n$. Then, $\Delta_{(i)i}=0$. 
Therefore the ideal $T_{i}= L_{i-1}+ \langle\Delta_{(i)i}\rangle= L_{i-1}$. We 
proceed in a similar way as Case 1. 
\bigskip

\section{Computation for $n=3$}
Let 
$X_{3}=
\begin{pmatrix}
0 & x_{12} & x_{13}\\
-x_{12} & 0 & x_{23}\\
-x_{13} & -x_{23} & 0\\
\end{pmatrix} $
 and $ Y_{3}=
\left[
\begin{array}{c}
y_{1} \\
y_{2} \\
y_{3} \\
\end{array}
\right]$. We write 
$$I_{1}(X_{3}Y_{3}) = \langle g_{13}, g_{23}, g_{33}\rangle,$$ 
such that 
\begin{eqnarray*}
g_{13} & = & x_{12}y_{2}+x_{13}y_{3}\\
g_{23} & = & -x_{12}y_{1}+x_{23}y_{3}\\
g_{33} & = & -x_{13}y_{1}-x_{23}y_{2}
\end{eqnarray*}
and $I_{3}=\langle g_{13}, g_{23}\rangle$; $J_{3}=\langle g_{33}\rangle$.
\medskip
 
We claim that $(I_{3}:J_{3})=\langle x_{12},y_{3}\rangle$. We first compute 
a Gr\"{o}bner basis of $I_{3}$. Let us fix the lexicographic monomial order induced by 
the ordering among the variables $y_{1}>y_{2}>y_{3}>x_{12}>x_{13}>x_{23}$ on $R$. 
Then $h=s(g_{13},g_{23})=x_{13}y_{3}y_{1}+x_{23}y_{3}y_{2}$. We have 
$\textrm{Lt}(h)=x_{13}y_{3}y_{1}$ and that it is not divisible by 
$\textrm{Lt}(g_{13})$ and  $\textrm{Lt}(g_{23})$. We therefore take 
the enlarged set $\{g_{13},g_{23},h\}$. It is clear that 
$\gcd(\textrm{LT}(g_{13}),\textrm{LT}(h))=1$, therefore we need 
to examine only $s(h,g_{23})$. Now 
$s(h,g_{23})=x_{13}x_{23}y_{3}^{2}+x_{12}x_{23}y_{3}y_{2}=x_{23}y_{3}(g_{13})\longrightarrow 0$; 
therefore the set $\{g_{13},g_{23},h\}$ forms a Gr\"{o}bner basis of $I_{3}$. 
We observe that $x_{12}g_{33}=x_{13}g_{23}-x_{23}g_{13}$ and $y_{3}g_{33}=-(y_{1}g_{13}+y_{2}g_{23})$. Therefore, $\langle x_{12},y_{3}\rangle \subset (I_{3}:J_{3}) $

Let $pg_{33}\in I_{3}$, and let $r$ be the remainder term upon division of $p$ by 
$x_{12},y_{3}$. We know that $\langle x_{12},y_{3}\rangle \subset (I_{3}:J_{3})$. 
Therefore, $rg_{33}\in I_{3}$. The set $\{g_{13},g_{23},h\}$ is a Gr\"{o}bner basis 
for $I_{3}$, therefore one of the following must hold: 
$x_{12}y_{2}\mid \textrm{Lt}(r)(x_{13}y_{1})$ 
or $x_{12}y_{1}\mid \textrm{Lt}(r)(x_{13}y_{1})$ or $x_{13}y_{1}y_{3}\mid \textrm{LT}(r)(x_{13}y_{1})$. This gives us $x_{12}\mid \textrm{Lt}(r)$ or $y_{3}\mid \textrm{Lt}(r)$, which 
leads to a contradiction if $r\neq 0$. Therefore $r=0$ and $p\in \langle x_{12},y_{3}\rangle$, 
and hence $(I_{3}:J_{3})=\langle x_{12},y_{3}\rangle$.
\medskip

Let $ L_{3}=I_{3}+J_{3} = \langle g_{13}, g_{23},g_{33} \rangle$. 
A minimal free resolution of $L_{3}$ is
$$0  \longrightarrow R^{2}\stackrel{d_{23}}\longrightarrow R^{3}\stackrel{d_{13}}{\longrightarrow} R\stackrel{d_{03}}{\longrightarrow}R/L_{3}\longrightarrow 0 $$
where $$d_{13}=\begin{pmatrix}

g_{13} & g_{23} & g_{33}\\
\end{pmatrix}, 
d_{23}=\begin{pmatrix}
x_{23} & y_{1}\\
-x_{13} & y_{2}\\
x_{12} & y_{3}\\
\end{pmatrix}$$
\bigskip
 
\section{Computation for $n=4$}

$X_{4}=
\begin{pmatrix}
0 & x_{12} & x_{13} & x_{14}\\
-x_{12} & 0 & x_{23} & x_{24}\\
-x_{13} & -x_{23} & 0 & x_{34}\\
-x_{14} & -x_{24} & -x_{34} & 0\\
\end{pmatrix} $
and 
$Y_{4}=\left[
\begin{array}{c}
y_{1} \\
y_{2} \\
y_{3} \\ 
y_{4} \\
\end{array}
\right]$
\medskip

\noindent By our notation we have,
\begin{eqnarray*}
g_{14} & = & x_{12}y_{2}+x_{13}y_{3}+x_{14}y_{4},\\ 
g_{24} & = & -x_{12}y_{1}+x_{23}y_{3}+x_{24}y_{4},\\ 
g_{34} & = & -x_{13}y_{1}-x_{23}y_{2}+x_{34}y_{4},\\ 
g_{44} & = & -x_{14}y_{1}-x_{24}y_{2}-x_{34}y_{3}
\end{eqnarray*}
and $I_{4}=\langle g_{14},g_{24},g_{34}\rangle $, $J_{4}=\langle g_{44}\rangle $, $L_{4}=I_{4}+J_{4}$, 
\medskip

We claim that, $C_{4}=(I_{4}:J_{4})=\langle g_{13},g_{23},g_{33}, y_{4}\rangle$. 
We first find a Gr\"{o}bner basis of $I_{4}$. Let us fix the lexicographic monomial 
order induced by $y_{1}>y_{2}>y_{3}>y_{4}>x_{12}>x_{13}>x_{14}>x_{23}>x_{24}>x_{34}$ 
on $R$. 
\medskip

Consider the $s$-polynomials: 
\begin{eqnarray*}
s(g_{14},g_{24}) & = & y_{1}y_{3}x_{13}+y_{1}y_{4}x_{14}+y_{2}y_{3}x_{23}+y_{2}y_{4}x_{23}\\
{} & = & -y_{3}g_{34}+y_{1}y_{4}x_{14}+y_{2}y_{4}x_{24}+y_{3}y_{4}x_{34}\\[2mm]
s(g_{24},g_{34}) & = & y_{2}x_{12}x_{23}+y_{3}x_{13}x_{23}-y_{4}x_{12}x_{34}+y_{4}x_{13}x_{24}\\
{} & = & x_{23}g_{14}+y_{4}x_{12}x_{34}-y_{4}x_{13}x_{24}+y_{4}x_{14}x_{23} 
\end{eqnarray*} 
 
\noindent We have $\gcd(\textrm{Lt}(g_{14}),\textrm{Lt}(g_{24})=1$, therefore $s(g_{14},g_{24})\longrightarrow 0$. Let us take 
$p_{1}=y_{1}y_{4}x_{14}+y_{2}y_{4}x_{24}+y_{3}y_{4}x_{34}$ and $p_{2}=y_{4}x_{12}x_{34}-y_{4}x_{13}x_{24}+y_{4}x_{14}x_{23}$ and consider the bigger set 
$\{g_{14},g_{24},g_{34},p_{1},p_{2}\}$. We now compute 
$$p_{3}=s(g_{14},p_{2})=y_{2}y_{4}x_{13}x_{24}-y_{2}y_{4}x_{14}x_{23}+y_{3}y_{4}x_{13}x_{34}+y_{4}^{2}x_{14}x_{34}.$$ 
It is evident that $\textrm{Lt}(p_{3})$ is not divisible by any element 
of the set 
$$\{\textrm{Lt}(g_{14}),\textrm{Lt}(g_{24}),\textrm{Lt}(g_{34}),\textrm{Lt}(p_{1}),\textrm{Lt}(p_{2})\}.$$ 
Therefore we add $p_{3}$ in the list and get the set $\mathcal{G}=\{g_{14},g_{24},g_{34},p_{1},p_{2},p_{3}\}$. It is now straightforward to check 
that every $s$ polynomial reduces to zero. Hence $\mathcal{G}$ is a a Gr\"{o}bner basis 
for the ideal $I_{4}$.
\medskip

We now compute a Gr\"{o}bner basis for the ideal 
$\langle g_{13},g_{23},g_{33}\rangle$. Consider the 
$s$-polynomials, 
\begin{eqnarray*}
s(g_{13},g_{23}) & = & x_{13}y_{3}y_{1}+x_{23}y_{3}y_{2} \quad = \quad -y_{3}g_{33}
\longrightarrow 0\\[2mm]
s( g_{23}, g_{33})& = & -y_{2}x_{12}x_{23}-y_{3}x_{13}x_{23} \quad = \quad -x_{23}g_{13}\longrightarrow 0.
\end{eqnarray*} 
Also, we have $\gcd(\textrm{Lt}(g_{13}),\textrm{Lt}(g_{33})=1$. Therefore, 
the set $\{g_{13},g_{23},g_{33}\}$ itself is a Gr\"{o}bner basis. Hence it 
follows easily that $\{g_{13},g_{23},g_{33},y_{4}\}$ is a Gr\"{o}bner basis 
for the ideal $\langle g_{13},g_{23},g_{33}, y_{4}\rangle$. 
\medskip

Using proposition \ref{useful} we obtain $\{g_{13},g_{23},g_{33},y_{4}\}\subset (I_{4}:J_{4})$, Let $pg_{44}\in I_{4}$ and assume that $r$ is the remainder upon division 
of $p$ by $\{g_{14},g_{24},g_{34},p_{1},p_{2}\}$. Suppose that $r\neq 0$. 
We have $rg_{44}\in I_{4}$. Moreover, $\textrm{Lt}(rg_{44})=\textrm{Lt}(r)x_{14}y_{1}$ 
is divisible by one of the leading terms $\textrm{Lt}(g_{14})=x_{12}y_{2}$, 
$\textrm{Lt}(g_{24})=x_{12}y_{1}$, $\textrm{Lt}(g_{34})=x_{13}y_{1}$, 
$\textrm{Lt}(p_{1})=y_{1}y_{4}x_{14}$, $\textrm{Lt}(p_{2})=y_{4}x_{12}x_{34}$, 
$\textrm{Lt}(p_{3})=y_{2}y_{4}x_{13}x_{24}$. If 
$\textrm{Lt}(rg_{44})$ is divisible by any one of the leading terms $\textrm{Lt}(g_{14})=x_{12}y_{2}$, $\textrm{Lt}(p_{1})=y_{1}y_{4}x_{14}$, $\textrm{Lt}(p_{2})=y_{4}x_{12}x_{34}$, 
$\textrm{Lt}(p_{3})=y_{2}y_{4}x_{13}x_{24}$, then we get a contradiction. If 
$\textrm{Lt}(g_{24})=x_{12}y_{1}\mid \textrm{LT}(rg_{44})$, then 
$x_{12}\mid \textrm{Lt}(r)$. Let $r=x_{12}m+l$. Therefore, 
$r.g_{44}=(x_{12}m+l)(-x_{14}y_{1}-x_{24}y_{2}-x_{34}y_{3})$ and 
after division we get 
$$q=(-x_{34}x_{12}y_{3}-x_{14}x_{23}y_{3}+x_{24}x_{13}y_{3})m+lg_{44}\in I_{4}.$$ 
We have $\textrm{Lt}(q)=x_{34}x_{12}y_{3}m$ and it must be divisible by one of the leading terms 
$\textrm{Lt}(g_{14})=x_{12}y_{2}$, $\textrm{Lt}(g_{24})=x_{12}y_{1}$, 
$\textrm{Lt}(g_{34})=x_{13}y_{1}$, $\textrm{Lt}(p_{1})=y_{1}y_{4}x_{14}$, 
$\textrm{Lt}(p_{2})=y_{4}x_{12}x_{34}$, $\textrm{Lt}(p_{3})=y_{2}y_{4}x_{13}x_{24}$. 
This implies that $\textrm{Lt}(r)$ must be divisible by one of the 
leading terms $\textrm{Lt}(g_{14})=x_{12}y_{2}$, $\textrm{Lt}(g_{24})=x_{12}y_{1}$, 
$\textrm{Lt}(g_{34})=x_{13}y_{1}$, $\textrm{Lt}(p_{1})=y_{1}y_{4}x_{14}$, 
$\textrm{Lt}(p_{2})=y_{4}x_{12}x_{34}$, $\textrm{Lt}(p_{3})=y_{2}y_{4}x_{13}x_{24}$, 
which is a contradiction. Similarly, if 
$\textrm{Lt}(g_{34})=x_{13}y_{1}\mid \textrm{Lt}(rg_{44})$, we get a 
contradiction. Therefore $r=0$ and our claim is proved.
\medskip

To find the resolution of $C_{4}$, we take the tensor product of the complexes:
$$0  \longrightarrow R^{2}\stackrel{d_{23}}\longrightarrow R^{3}\stackrel{d_{13}}{\longrightarrow} R\stackrel{d_{03}}{\longrightarrow}R/L_{3}\longrightarrow 0 $$
and
$$0  \longrightarrow R\stackrel{y_{4}}\longrightarrow R{\longrightarrow}R/y_{4}R\longrightarrow 0 $$
and obtain a resolution of $C_{4}$ as
$$0  \longrightarrow R^{2}\stackrel{\eta_{34}}\longrightarrow R^{5}\stackrel{\eta_{24}}{\longrightarrow} R^{4}\stackrel{\eta_{14}}{\longrightarrow}R\stackrel{\eta_{04}}{\longrightarrow} R/C_{4}\longrightarrow 0 $$
where
\begin{itemize}
\item $\eta_{14}=[d_{13}|y_{4}]=[g_{13},g_{23},g_{33}, y_{4}]$,
\medskip

\item $\eta_{24} = 
\left[
\begin{array}{c|c}
d_{23} & -y_{4}I_{3}\\ \hline
0 & d_{13}
\end{array}\right] 
= 
\begin{pmatrix}
x_{23} & y_{1} & -y_{4} & 0 & 0\\
-x_{13} & y_{2} & 0 & -y_{4} & 0\\
x_{12} & y_{3} & 0 & 0 & -y_{4}\\
0 & 0 & g_{13} & g_{23} & g_{33} \\
\end{pmatrix}$,
\medskip

\item $\eta_{34}=
\begin{pmatrix}
y_{4} & 0 \\
0 & y_{4} \\
x_{23} & y_{1} \\
-x_{13} & y_{2} \\
x_{12} & y_{3} \\
\end{pmatrix}$.
\end{itemize}
\medskip

\noindent Using the mapping cone between these complexes we get
$$\xymatrix{
0 \ar[r] &R^{2}\ar[d]^{\xi{34}}\ar[r]^{\eta{34}}& R^{5}\ar[r]^{\eta_{24}}\ar[d]^{\xi_{24}}&
R^{4}\ar[r]^{\eta_{14}}\ar[d]^{\xi_{14}}&R\ar[r]^{\eta_{04}}\ar[d]^{\xi_{04}=g_{44}} &R/C_{4}\ar[d]\ar[r]&0 \\
0 \ar[r] &R\ar[r]_{\psi_{34}} &R^{3}\ar[r]_{\psi_{24}}&R^{3}\ar[r]_{\psi_{14}}&R\ar[r]_{\psi_{04}}&R/I_{4}\ar[r]& 0
} 
$$
where
\begin{itemize}
\item $\psi_{14}=[g_{14},g_{24},g_{34}]$, \,\,$\psi_{24}=
\begin{pmatrix}
g_{24} & 0 & -g_{34}\\
-g_{14} & -g_{34} & 0\\
0 & -g_{24} & g_{14} \\
\end{pmatrix}$,
\medskip

\item  $\psi_{34}=\left[
\begin{array}{c}
g_{34} \\
g_{14} \\
g_{24} \\ 
\end{array}
\right]$
\medskip

\item $\xi_{04}=[g_{44}]$,
\medskip

\item $\xi_{14}=
\begin{pmatrix}
g_{44}+x_{14}y_{1} & x_{24}y_{1} & x_{34}y_{1} & -y_{1}\\
x_{14}y_{2} & g_{44}+x_{24}y_{2} & x_{34}y_{2} & -y_{2}\\
x_{14}y_{3} & x_{24}y_{3} & g_{44}+x_{34}y_{3} & -y_{3}\\
\end{pmatrix}$,
\medskip

\item $\xi_{24}=
\begin{pmatrix}
-x_{34} & 0 & y_{2} & -y_{1} & 0\\
-x_{14} & 0 & 0 & y_{3} & -y_{2}\\
-x_{24} & 0 & -y_{3} & 0 & y_{1}\\
 \end{pmatrix}$,
\medskip

\item $\xi_{24}=
\begin{pmatrix}
 -1 & 0\\
\end{pmatrix}$.
\end{itemize}
\medskip

\noindent Therefore a non-minimal resolution of $L_{4}$ is
$$ 0  \longrightarrow R^{2}\stackrel{\widetilde{d_{44}}}\longrightarrow R^{6}\stackrel{\widetilde{d_{34}}}{\longrightarrow} R^{7}\stackrel{\widetilde{d_{24}}}{\longrightarrow}R^{4}\stackrel{\widetilde{d_{14}}}{\longrightarrow}R\stackrel{\widetilde{d_{04}}}{\longrightarrow}R/L_{4}\longrightarrow 0 , $$
where $\widetilde{d_{14}}=[g_{14},g_{24},g_{34},g_{44}]$, 
$\quad \widetilde{d_{24}}=
\left[
\begin{array}{c|c}
\xi_{14} & \psi_{24} \\ \hline \\
-\eta_{14} & 0
\end{array}\right]$, 

$\quad \widetilde{d_{34}}=
\left[
\begin{array}{c|c}
-\eta_{24}  & 0\\ \hline
\xi_{24}  & \psi_{34}
\end{array}\right]$, $\quad \widetilde{d_{44}}=\left[
\begin{array}{c}
-\eta_{34}\\ \hline
\xi_{34}\\
\end{array}
\right]$. 
\medskip

\noindent Therefore a minimal free resolution of $L_{4}=I_{1}(X_{4}Y_{4})$ is
$$ 0  \longrightarrow R\stackrel{d_{44}}\longrightarrow R^{5}\stackrel{d_{34}}{\longrightarrow} R^{7}\stackrel{d_{24}}{\longrightarrow}R^{4}\stackrel{d_{14}}{\longrightarrow}R\stackrel{d_{04}}{\longrightarrow}R/L_{4}\longrightarrow 0 $$
where $d_{14}=\widetilde{d_{14}}$, $\quad d_{24}=\widetilde{d_{24}}$, 
$\quad d_{34}=
\left[
\begin{array}{c}
-\eta_{24}\\ \hline
\xi_{24} 
\end{array}\right]$, $\quad d_{44}=\left[
\begin{array}{c}
0\\
y_{4}\\
y_{1}\\
y_{2}\\
y_{3}
\end{array}
\right]$.
\bigskip

\bibliographystyle{amsalpha}

\end{document}